\documentclass[11pt,draft,a4paper]{amsart}
\usepackage{amsfonts,amsmath,amssymb,amscd}

\newtheorem{theorem}{Theorem}

\begin{document}
\title{Cellular cochain algebras and torus actions}
\address{Faculty of Mathematics and Mechanics, Moscow State University}
\author{Ilia~V.~Baskakov}
\author{Victor~M.~Buchstaber}
\author{Taras~E.~Panov}
\thanks{The work is supported by the Russian Foundation for Basic
Research, grant no.~02-01-00659, and by the Leading Scientific School
support, grant no.~2185.2003.1.}
\email{tpanov@mech.math.msu.su}
\maketitle

We prove that the \emph{integral} cohomology algebra of the moment-angle complex
$\mathcal Z_K$~\cite{bu-pa02} is isomorphic to the $\operatorname{Tor}$-algebra of the face ring 
of simplicial complex~$K$. The proof relies upon the construction of a cellular approximation
of the diagonal map
$\Delta\colon\mathcal Z_K\to\mathcal Z_K\times\mathcal Z_K$.
Cellular cochains do not admit a functorial associative multiplication because a proper cellular
diagonal approximation does not exist in general. The construction of moment-angle complexes is a functor 
from the category of simplicial complexes to the category of spaces with torus action. We show 
that in this special case the proposed cellular approximation of the diagonal is associative and
functorial with respect to those maps of moment-angle complexes which are induced by simplicial maps.

The \emph{face ring} of a complex $K$ on the vertex set $[m]=\{1,\dots,m\}$ is the graded quotient ring
$\mathbb Z[K]=\mathbb Z[v_1,\dots,v_m]/(v_\omega\colon\omega\notin K)$
with $\deg v_i=2$ and $v_\omega=v_{i_1}\dotsb v_{i_k}$, where
$\omega=\{i_1,\dots,i_k\}\subseteq[m]$. Let
$BT^m$ be the classifying space for the $m$-dimensional torus, endowed with the standard cell decomposition.
Consider a cellular subcomplex
$\text{\it DJ\/}(K):=\bigcup_{\sigma\in K}BT^\sigma\subseteq BT^m$,
where
$BT^\sigma=\{x=(x_1,\dots,x_m)\in BT^m:x_i=pt\text{ for }i\notin\sigma\}$.
Using the cellular decomposition we establish a ring isomorphism
$H^*(\text{\it DJ\/}(K))\cong\mathbb Z[K]$ (see \cite[Lemma~2.8]{bu-pa00}).
Let $D^2\subset\mathbb C$ be the unit disc and set
$B_\omega:=\{(z_1,\dots,z_m)\in(D^2)^m:|z_j|=1\text{ for }j\notin\omega\}$.
The \emph{moment-angle complex} is the $T^m$-invariant subspace
$\mathcal Z_K:=\bigcup_{\sigma\in K}B_\sigma\subseteq(D^2)^m$.
As it is shown in~\cite[Ch.~6]{bu-pa02}, the spaces
$\text{\it DJ\/}(K)$ and $\mathcal Z_K$ are homotopy equivalent to the spaces introduced 
in~\cite{da-ja91}, which justifies our notation. Complexes 
$\mathcal Z_K$ provide an important class of torus actions. The space $\mathcal Z_K$ is 
the homotopy fibre of the inclusion $\text{\it DJ\/}(K)\hookrightarrow BT^m$; it appears also as
the level surface of the moment map used in the construction of \emph{toric varieties} via symplectic reduction; it is also homotopy equivalent to the complement of the \emph{coordinate subspace arrangement} determined by~$K$, see~\cite[\S\,8.2]{bu-pa02}.

\begin{theorem}
There is a functorial in $K$ isomorphism of algebras
$$
H^*(\mathcal Z_K;\mathbb Z)\cong
H\bigl[\Lambda[u_1,\dots,u_m]\otimes\mathbb Z[K],d\bigr]
\cong\operatorname{Tor}_{\mathbb Z[v_1,\dots,v_m]}\bigl(\mathbb Z[K],\mathbb Z\bigr),
$$
where the algebra in the middle is the cohomology of the differential graded algebra with
$\deg u_i=1$, $\deg v_i=2$, $du_i=v_i$, $dv_i=0$.
\end{theorem}

In the case of rational coefficients this theorem was proved in~\cite{bu-pa99} using spectral sequences techniques
(see also \cite[Th.~7.6, Probl~8.14]{bu-pa02}). Our new proof uses a construction of cellular cochain algebra.
Another proof of Theorem~1 follows from a recent independent work of M.~Franz~\cite[Th.~1.2]{fran??}.

\begin{proof}[Proof of Theorem 1]
We prove just the first isomorphism, since the second is a standard consequence 
of the Koszul resolution (see~\cite{bu-pa02} for details). Introduce an extra grading by 
setting $\operatorname{bideg}u_i=(-1,2)$, $\deg v_i=(0,2)$, and consider a quotient algebra
$$
R^*(K):=\Lambda[u_1,\dots,u_m]\otimes\mathbb Z[K]\big/
(v_i^2=u_iv_i=0,\ i=1,\dots,m).
$$
Let $\varrho\colon\Lambda[u_1,\dots,u_m]\otimes\mathbb Z[K]\to R^*(K)$ be the canonical
projection. Using the finite additive basis in $R^*(K)$ consisting of monomials of the type
$u_\omega v_\sigma$, where $\omega\subseteq[m]$, $\sigma\in K$ and
$\omega\cap\sigma=\varnothing$, we define an additive inclusion
$\iota\colon R^*(K)\to\Lambda[u_1,\dots,u_m]\otimes\mathbb Z[K]$
which satisfies $\varrho\cdot\iota=\operatorname{id}$. We claim that 
$\varrho$ induces an isomorphism in the cohomology. To see this we construct 
a cochain homotopy operator $s$ such that $ds+sd=\operatorname{id}-\iota\cdot\varrho$. If 
$K=\varDelta^{m-1}$ (the full simplex) the algebra
$\Lambda[u_1,\dots,u_m]\otimes\mathbb Z[K]$ is
$E=E_m=\Lambda[u_1,\dots,u_m]\otimes\mathbb Z[v_1,\dots,v_m]$,
and $R^*(\varDelta^{m-1})$ is isomorphic to the algebra
$R^*(\varDelta^0)^{\otimes m}$, where
$R^*(\varDelta^0)=
\Lambda[u]\otimes\mathbb Z[v]\bigr/(v^2=uv=0)$.
A direct calculation shows that for $m=1$ the map
$s_1\colon E^{0,*}=\mathbb Z[v]\to E^{-1,*}$ given by
$s_1(a_0+a_1v+\dots+a_jv^j)=(a_2v+a_3v^2+\dots+a_jv^{j-1})u$
is the required cochain homotopy. Then we may assume by induction that for $m=k-1$ 
there is a cochain homotopy operator $s_{k-1}\colon E_{k-1}\to E_{k-1}$. Since $E_k=E_{k-1}\otimes E_1$,
$\varrho_k=\varrho_{k-1}\otimes\varrho_1$ and
$\iota_k=\iota_{k-1}\otimes\iota_1$, the map
$s_k=s_{k-1}\otimes\operatorname{id}+\iota_{k-1}\varrho_{k-1}\otimes s_1$
is a cochain homotopy between $\operatorname{id}$ and
$\iota_k\varrho_k$, which finishes the proof for $K=\varDelta^{m-1}$. In the case of arbitrary complex
$K$ the algebras
$\Lambda[u_1,\dots,u_m]\otimes\mathbb Z[K]$ and
$R^*(K)$ are obtained from $E_m$ and $R^*(\varDelta^{m-1})$ respectively by factoring out a monomial ideal. 
This factorisation does not affect the properties of the operator $s$, which therefore establishes the required
cochain homotopy.

Consider a cellular subdivision of the polydisc~$(D^2)^m$ with each $D^2$ subdivided
into the cells 1, $T$ and $D$ of dimensions $0$, $1$ and~$2$ respectively. Therefore,
the cells of $(D^2)^m$ are encoded by words $\mathcal T\in\{D,T,1\}^m$. Assign to each
pair of subsets $\sigma,\omega\subseteq[m]$, $\sigma\cap\omega=\varnothing$
the word $\mathcal T(\sigma,\omega)$ which has the letter $D$ on the positions numbered by $\sigma$
and letter $T$ on the positions with numbers from~$\omega$. By the construction, the word
$\mathcal T(\sigma,\omega)$ corresponds to a cell from $\mathcal Z_K\subset(D^2)^m$
if and only if $\sigma\in K$. The cellular cochain complex
$C^*(\mathcal Z_K)$ has an additive basis consisting of the cochains of the type
$\mathcal T(\sigma,\omega)^*$. It follows that the dimensions of the graded
components of $C^*(\mathcal Z_K)$ and $R^*(K)$ coincide. Moreover, the map
$g\colon R^*(K)\to C^*(\mathcal Z_K)$,
$u_\omega v_\sigma\mapsto\mathcal T(\sigma,\omega)^*$,
is an isomorphism of differential graded modules and induces an additive isomorphism
$H[R^*(K)]\cong H^*(\mathcal Z_K)$. To finish the proof we construct a cellular
approximation
$\widetilde\Delta_K$ of the diagonal map
$\Delta\colon\mathcal Z_K\to\mathcal Z_K\times\mathcal Z_K$
such that the multiplication induced on the cellular cochains coincides with the 
multiplication in $R^*(K)$ via the isomorphism~$g$. For $K=\varDelta^0$ we have
$\mathcal Z_K=D^2$. Set $z=\rho e^{i\varphi}\in D^2$ and define a map
$\widetilde\Delta\colon D^2\to D^2\times D^2$ by
$$
\widetilde\Delta(z)=\begin{cases}
(1+\rho(e^{2i\varphi}-1),1)&\text{for }\varphi\in[0,\pi],\\
(1,1+\rho(e^{2i\varphi}-1))&\text{for }\varphi\in[\pi,2\pi).
\end{cases}
$$
This is a cellular map sending $\partial D^2$ to
$\partial D^2\times\partial D^2$ and homotopic to $\Delta$ in the class of such maps.
For the corresponding multiplication in cellular cochains there is a multiplicative
isomorphism $R^*(\varDelta^0)\to C^*(D^2)$.
Passing to the case $K=\varDelta^{m-1}$, we see that the map
$\widetilde\Delta\colon(D^2)^m\to(D^2)^m\times(D^2)^m$ gives rise to a 
multiplicative isomorphism
\begin{multline*}
f\colon R^*(\varDelta^{m-1})=
\Lambda[u_1,\dots,u_m]\otimes\mathbb Z[v_1,\dots,v_m]\bigl/(v_i^2=u_iv_i=0)\\
\longrightarrow C^*\bigl((D^2)^m\bigr).
\end{multline*}
As it follows from the construction of $\mathcal Z_K$, the restriction of the
map $\widetilde\Delta$ to $\mathcal Z_K$ is a cellular map
$\widetilde\Delta_K\colon\mathcal Z_K\to\mathcal Z_K\times\mathcal Z_K$.
Therefore, there is a multiplicative map
$q\colon C^*((D^2)^m)\to C^*(\mathcal Z_K)$. Consider the commutative diagram
$$
\begin{CD}
R^*(\varDelta^{m-1})@>f>>C^*\bigl((D^2)^m\bigr)\\
@VpVV@VVqV\\
R^*(K)@>g>>C^*(\mathcal Z_K)
\end{CD}
\qquad.
$$
Using the fact that $p$, $f$ and $q$ are ring homomorphisms, and
$p$ is an epimorphism, we deduce that $g$ is also a ring homomorphism. Thus,
$g$ is a multiplicative isomorphism.
\end{proof}

\end{document}